\documentclass[12pt]{amsart}
\usepackage{amsfonts,latexsym,amsthm,amssymb}
\usepackage[all]{xy}

\newtheorem{theorem}{Theorem}[section]
\newtheorem{lemma}[theorem]{Lemma}
\newtheorem{corol}[theorem]{Corollary}
\newtheorem{prop}[theorem]{Proposition}

{\theoremstyle{definition} }
{\theoremstyle{remark} 
}

\newcommand{\Abb}{{\mathbb{A}}}
\newcommand{\Cbb}{{\mathbb{C}}}

\newcommand{\Lbb}{{\mathbb{L}}}
\newcommand{\Pbb}{{\mathbb{P}}}
\newcommand{\Qbb}{{\mathbb{Q}}}
\newcommand{\Zbb}{{\mathbb{Z}}}

\newcommand{\Til}[1]{{\widetilde{#1}}}

\DeclareMathOperator{\Spec}{Spec}

\newcommand{\one}{1\hskip-3.5pt1}
\newcommand{\csm}{c_{\text{\rm SM}}}

\title{Chern classes of birational varieties}
\author{Paolo Aluffi}
\address{Florida State University, Tallahassee, Florida}
\email{aluffi@math.fsu.edu}

\begin{document}

\begin{abstract}
Let $\varphi: V\dashrightarrow W$ be a birational map between smooth
algebraic varieties which does not change the canonical class (in the
sense of Batyrev, \cite{MR2000i:14059}). We prove that the total
homology Chern classes of $V$ and $W$ are push-forwards of the same
class from a resolution of indeterminacies of $\varphi$.

For example, it follows that the push-forward of the total Chern class
of a crepant resolution of a singular variety is independent of the
resolution.
\end{abstract}

\maketitle

\section{Introduction and statement of the result}\label{intro}

There is a strong motivic feel to the theory of
Chern-Schwartz-MacPherson classes, although this does not seem to have
yet congealed into a precise statement in the literature. In this note
we give an instance of a result exploiting this heuristic observation.
Our statement is an analog of Victor Batyrev's well-known theorem
showing that the Betti numbers of birational manifolds `in the same
$K$-class' are equal, see \cite{MR2000i:14059}. Batyrev's theorem has
been generalized to many other numerical invariants, such as
the Hodge numbers or certain Chern numbers. Our statement is not
numerical; it refers to the total Chern class of the tangent bundle, a
class in the Chow group of the variety.

\begin{theorem}\label{main}
Let $\varphi: V \dashrightarrow W$ be a birational morphism of
nonsingular algebraic varieties over an algebraically closed field of
characteristic~0. Assume that there is a resolution of
indeterminacies of $\varphi$, $Z$:
$$\xymatrix{
& Z \ar[ld]_v \ar[rd]^w \\
V \ar@{-->}[rr] & & W
}$$
such that $v$, $w$ are proper and birational, and the Jacobian ideals
of $v$ and $w$ coincide. 

Then there exists a class $C\in (A_*Z)_\Qbb$ such that
$$c(TV)\cap [V] = v_*(C)\quad\text{and}\quad c(TW)\cap [W]= w_*(C)$$
in $(A_*V)_\Qbb$, $(A_*W)_\Qbb$ respectively.
\end{theorem}

The condition on the Jacobians implies that the pull-backs $v^* K_V$
and $w^* K_W$ of the respective canonical classes coincide. In fact,
as Lawrence Ein kindly pointed out to me, the condition {\em is
equivalent\/} to requiring that $v^* K_V \equiv_{\text{Num}} w^* K_W$,
an easy consequence of the `Kodaira lemma' (see for example 2.19 in
\cite{MR94f:14013}). 

Theorem~\ref{main} has a number of immediate consequences: for instance
the equality of Euler characteristics of $V$ and $W$, or more generally
the equality of Chern numbers
$$c_1(V)^i\cdot c_{n-i}(V)=c_1(W)^i\cdot c_{n-i}(W)$$
for all $i\ge 0$ (where $n=\dim V=\dim W$) in the stated hypotheses,
when $V$ and $W$ are complete. Identities of this type are not new;
for example, Anatoly Libgober and John Wood have proved
(\cite{MR91g:32039}) that $c_1\cdot c_{n-1}$ is determined by the
Hodge numbers, and these agree for $V$ and $W$, as mentioned above. 
Our argument is very direct, and yields the equality for all $i$;
suprisingly, this simple equality does not seem to be explicitly in
the literature for $i\ge 2$, although it is certainly implied by 
recent and powerful results on elliptic genera (\cite{MR1953295},
\cite{MR2003j:14015}).  Also:

\begin{corol}\label{crep}
Let $\alpha: Y \to X$ be a crepant resolution. Then the class
$$\alpha_*(c(TY)\cap [Y])$$
in $(A_*X)_\Qbb$ is independent of $Y$.
\end{corol}

Again, this is immediate from Theorem~\ref{main}; it also follows from
the more sophisticated technology in Lev Borisov and Libgober's work. 
It would be interesting to extend this notion of `Chern class' to
arbitrary singular varieties, and to compare it to other known
intrinsic definitions such as the Chern-Fulton class
(\cite{MR85k:14004}, Example~4.2.6), or the Chern-Schwartz-MacPherson
class (\cite{MR50:13587}).

Similarly to analogous numerical statements, Theorem~\ref{main} is an
easy consequence of a basic formula in motivic integration, to wit
Proposition~6.3.2 in \cite{MR99k:14002}. This formula is usually
obtained as a corollary of Kontsevich's change of variable formula; we
refer the reader to loc.~cit.~or to the excellent surveys
\cite{MR1905328}, \cite{MR2003k:14010} for a discussion in the context
of motivic integration. 

In order to provide a self-contained treatment of the material, we
take the liberty of providing an alternative proof of this `motivic'
statement, using the factorization theorem of \cite{MR2003c:14016};
our argument has the (very small) advantage of proving the needed
formula without requiring completions. This is done in
\S\ref{motiv}. The idea of using the factorization theorem as an
alternative to motivic integration is of course not new; for example,
Willem Veys uses it to prove a similar (in fact, stronger) statement,
in \S2 of \cite{MR2030094}. Also, the factorization theorem is a
key tool in the work of Borisov and Libgober on elliptic genera. In
fact, Borisov has pointed out to me that Theorem~\ref{main} can also
be derived from Theorem~3.5 in \cite{math.AG/0206241}, and that this
work may be used to extend the class given in Corollary~\ref{crep} to
all varieties with log-terminal singularities.

Regardless of how the basic formula is established, however, our main
observation is that this formula alone---and in fact even just the
Euler characteristic version of the basic formula, which goes back to
\cite{MR93g:11118}---implies an analog at the level of
Chern-Schwartz-MacPherson classes, and that this immediately implies
Theorem~\ref{main}. The upgrade of the basic formula to
Chern-Schwartz-MacPherson classes is Theorem~\ref{strmain};
Theorem~\ref{main} follows easily from this result, as shown in
\S\ref{CSM}.

Theorem~\ref{strmain} must thus have been known to the experts very
early (and Fran\c cois Loeser confirms this); the fact that it has
nice applications, such as the invariance of $c_1^i\cdot c_{n-i}$
under $K$-equivalence, seems to have been inexplicably overlooked thus
far. It may be argued that Theorem~\ref{main}, and hence these
straightforward consequences, should have been immediately noticed
after the diffusion of \cite{MR93g:11118}; in particular, long before
the discovery of the factorization theorem of \cite{MR2003c:14016},
which has fueled the more recent and sophisticated work in this
direction. 

In any case, we believe that Theorem~\ref{strmain} is of independent
interest in the context of Chern-Schwartz-MacPherson classes, and that
it deserves to be out in the open.

I thank Prof.~Kenji Matsuki for promptly answering many queries
concerning \cite{MR2003c:14016}. 
Thanks are also due to the Max-Planck-Institut f\"ur Mathematik in 
Bonn, where part of this work was done.

\section{The basic formula}\label{motiv}

We work over an algebraically closed field $k$. We will assume $k$ has
characteristic~0, as this is at present needed for resolution of
singularities (needed in the factorization theorem of
\cite{MR2003c:14016}) and the theory of Chern-Schwartz-MacPherson
classes used in \S\ref{CSM}.

If $X$ is a variety, $[X]$ will denote its class in the Grothendieck
ring of algebraic varieties over $k$; that is, the abelian group
generated by isomorphism classes of varieties modulo the relation
$[X-Y]=[X]-[Y]$ whenever $Y$ is a closed subvariety of $X$, with
multiplication defined by product over~$k$. The unit element is
$1=[\Spec k]$. It is costumary to denote by $\Lbb$ the `Tate motive'
$[\Abb^1]$.

The basic formula is a relation in this ring, depending on the datum
of a `modification' of nonsingular varieties. It is in fact
convenient to localize the ring at elements
$[\Pbb^\mu]=\frac{\Lbb^{\mu+1}-1}{\Lbb-1}$, and this will be done
without further mention in what follows.

{\bf Notation:} if $\{D_j\}_{j\in J}$ is a set of irreducible divisors
in a variety and $I\subset J$, then $D_I^\circ$ will denote the
complement of $\cup_{j\in J\setminus I}D_j$ in $\cap_{i\in I}D_i$.

Let $v: Z \to V$ be a proper birational morphism of nonsingular
varieties, such that the exceptional divisor of $v$ has normal
crossings, with nonsingular irreducible components $E_j$, $j\in J$.

Assume that the Jacobian ideal of $v$ is principal, with divisor
$\sum_j \mu_j E_j$, so that $v^*K_V=K_Z+\sum_j \mu_j E_j$. The
following formula is (essentially) a particular case of Proposition~6.3.2 
in \cite{MR99k:14002}.

\begin{theorem}\label{basic}
Let $U$ be a closed subvariety of $V$. Then
$$[U]=\sum_{I\subset J}\frac{[E_I^\circ\cap v^{-1}(U)]}{\prod_{i\in I}
[\Pbb^{\mu_i}]}\quad.$$
\end{theorem}

As mentioned in \S\ref{intro}, we provide a motivic-integration-free
proof of Theorem~\ref{basic} in this section, which the hurried reader
may safely skip.

The factorization theorem of \cite{MR2003c:14016} will reduce the
proof of Theorem~\ref{basic} to the case in which $v$ is a blow-up at
a smooth center; this case will be dealt with by analyzing the
situation in the exceptional divisor, and this in turn will be reduced
to the analysis of a single fiber, a projective space. The real reason
behind Theorem~\ref{basic} can then be traced to the following simple
lemma, whose statement requires no localization.

\begin{lemma}\label{simplex}
Let $d\ge k>0$ and $\mu_i\ge 0$, $i=1,\dots, k$, be integers, and let
$K=\{1,\dots,k\}$. Let $H_1,\dots,H_k$ be linearly independent
hyperplanes in $\Pbb^{d-1}$. 
Then
$$\sum_{I\subset K} \left([H_I^\circ]\prod_{i\in K\setminus 
I}[\Pbb^{\mu_i}]\right)=[\Pbb^{\sum_{j\in K}\mu_j+d-1}]\quad.$$
\end{lemma}

\begin{proof}
With evident coordinates, the subset $H_I^\circ$ consists of
homogeneous $d$-tuples in which $|I|$ fixed components are~0 and
$k-|I|$ fixed components are nonzero. It follows that
$$[H_I^\circ]=\left\{
\aligned
\frac{(\Lbb-1)^{k-|I|}\Lbb^{d-k}}{\Lbb-1} &\quad |I|\ne k\\
\frac{\Lbb^{d-k}-1}{\Lbb-1}  \qquad &\quad |I|=k
\endaligned
\right.\quad,$$
giving
$$[H_I^\circ]\prod_{i\in K\setminus I}[\Pbb^{\mu_i}]=\left\{
\aligned
\frac{\Lbb^{d-k}}{\Lbb-1} \prod_{i\in K\setminus I}(\Lbb^{\mu_i+1}-1) 
&\quad I\ne K\\
\frac{\Lbb^{d-k}-1}{\Lbb-1} - \frac1{\Lbb-1} \quad &\quad I=K
\endaligned
\right.\quad.$$
Hence the left-hand-side in the stated equality is
$$\frac{\Lbb^{d-k}}{\Lbb-1}\sum_{I\subset K} \prod_{i\in K\setminus I} 
(\Lbb^{(\mu_i+1)}-1) - \frac 1{\Lbb-1}
= \frac{\Lbb^{d-k}}{\Lbb-1}\prod_{j\in K} (1+ (\Lbb^{\mu_j+1}-1)) - 
\frac 1{\Lbb-1}$$
which immediately yields the right-hand-side.
\end{proof}

The form in which Lemma~\ref{simplex} will be used is the following.
\begin{corol}\label{simplexcor}
With the same notation as in Lemma~\ref{simplex}, set 
$\mu_0=\sum_{j\in K}\mu_j+d-1$; then
$$\sum_{I\subset K} \frac{[H_I^\circ]}
{\prod_{i\in (\{0\}\cup I)}[\Pbb^{\mu_i}]}
=\frac 1{\prod_{j\in K}[\Pbb^{\mu_j}]}\quad.$$
\end{corol} 

This formula will be used to control the situation across blow-ups. 
Let $X$ be a nonsingular algebraic variety, and assume given a set
$\{E_j\}_{j\in J}$ of irreducible nonsingular divisors in $X$. Let
$S\subset X$ be a nonsingular subvariety of codimension~$d\ge 1$,
intersecting $\cup E_j$ with normal crossings: that is, at each point
$s\in S$ there is an analytic system of parameters $x_1,\dots,x_n$ for
$X$ at $s$ such that $S$ is given by $x_1=\dots=x_d=0$, and $\cup E_j$
is given by a monomial in the $x_i$'s. 

Let $\pi:Y\to X$ be the blow-up of $X$ along $S$; denote by $\Til E_0$
the exceptional divisor $\pi^{-1}(S)$, and let $\Til E_j$ be the
proper transform of $E_j$. The following lemma is a straightforward
computation, which we leave to the reader.

\begin{lemma}\label{standard}
---The divisor $\Til E_0\cup(\cup_{j\in J} \Til E_j)$ has normal
crossings at all points of $\Til E_0$;

---If $E_j$ does not contain $S$ and $s\in S\cap E_j$, then $\Til
E_j$ contains the fiber $\pi^{-1}(s)$;

---Let $K$ denotes the set of indices $j\in J$ such that $E_j$
contains $S$. Then for all $s\in S$ and $j\in K$, the intersections of
$\Til E_j$ with $\pi^{-1}(s)\cong \Pbb^{d-1}$ consist of linearly
independent hyperplanes in $\Pbb^{d-1}$.
\end{lemma}

Now assign to each $E_j$ and $\Til E_j$, $j\in J$ a nonnegative integer
$\mu_j$, and assign to the exceptional divisor $\Til E_0$ the combination
$$\mu_0:=\sum_{E_j\supset Z}\mu_j+d-1\quad.$$
Note that, with this choice,
$$(\sum_{j\in (\{0\}\cup J)}\mu_j \Til E_j) - \pi^* (\sum_{j\in J}
\mu_j E_j)$$
is the divisor of the Jacobian ideal of $\pi$.

For any algebraic subset $T\subset X$, let
$$\chi_X(T):=\sum_{I\subset J}\frac{[E_I^\circ\cap T]}
{\prod_{i\in I}[\Pbb^{\mu_i}]}$$
and
$$\chi_Y(\pi^{-1}(T)):=\sum_{I\subset (\{0\}\cup J)}\frac
{[\Til E_I^\circ\cap \pi^{-1}(T)]}{\prod_{i\in I}[\Pbb^{\mu_i}]}
\quad.$$

\begin{prop}\label{pre}
With notation as above, $\chi_Y(\pi^{-1}(T))=\chi_X(T)$.
\end{prop}

\begin{proof}
The contributions to $\chi_X(T)$ coming from the complement of $S$
equal the contributions to $\chi_Y(\pi^{-1}(T))$ coming from the
complement of $\Til E_0$. Thus it suffices to prove that 
$$\sum_{I\subset J}\frac{[E_I^\circ\cap T\cap S]}{\prod_{i\in I}
[\Pbb^\mu_i]}$$ 
equals the contribution to $\chi_Y(\pi^{-1}(T))$ supported
within $\Til E_0$, which consists of
$$\sum_{I\subset J}\frac{[\Til E_{\{0\}\cup I}^\circ\cap
  \pi^{-1}(T)]}{[\Pbb^{\mu_0}]\prod_{i\in I}[\Pbb^\mu_i]}\quad.$$
Now every point of $T\cap S$ belongs to exactly one $E_K^\circ$.
For all $K\subset J$, $\Til E_{\{0\}\cup K}^\circ\cap
  \pi^{-1}(T)$ fibers onto $E_K^\circ\cap T\cap S$ with fibers equal
to the union over all $I\subset K$ of
$$\Til E_I^\circ\cap \pi^{-1}(x) \cong \Til E_I^\circ\cap \Pbb^{d-1}$$
for $x\in E_K^\circ\cap T\cap S$ (identifying the fiber of $\Til E_0$
with $\Pbb^{d-1}$). By the multiplicativity of products in the
Grothendieck ring, it suffices to show that for all $K\subset J$
$$\sum_{I\subset K}\frac{[\Til E_I^\circ \cap \pi^{-1}(x)]}
{[\Pbb^{\mu_0}]\prod_{i\in I}[\Pbb^{\mu_i}]} = \frac 1{\prod_{j\in K}
  [\Pbb^{\mu_j}]}\quad,$$
for all $x\in T\cap S\cap E_K^\circ$. By Lemma~\ref{standard}, indices
in $K$ corresponding to divisors which do not contain $S$ only
contribute common multiplicative factors to the denominators in both
sides of this equality; so we may assume that, for all $j$ in $K$,
$S\subset E_j$.

Again by Lemma~\ref{standard}, the intersections $\pi^{-1}(x)\cap \Til
E_j$ consist of linearly independent hyperplanes in $\pi^{-1}(x)\cong
\Pbb^{d-1}$; with the given $\mu_0$, the sought formula is then
precisely the statement of Corollary~\ref{simplexcor}.
\end{proof}

We are now ready to prove Theorem~\ref{basic}. 

\begin{proof}
By the factorization theorem of \cite{MR2003c:14016}, $v$ may be
decomposed as a sequence:
$$\xymatrix{
Z=V_0 \ar@{-->}[r] & V_1 \ar@{-->}[r] & V_2 \ar@{-->}[r] & \dots
\ar@{-->}[r] & V_M=V
}$$
where each map is either a blow-up or a blow-down at a smooth center.
More precisely, for each $i$ one has one of the possibilities
$$\xymatrix{
V_i \ar[rr] \ar[rd] & & V_{i+1} \ar[ld] \\
& V
}\qquad
\xymatrix{
V_{i+1} \ar[rr] \ar[rd] & & V_i \ar[ld] \\
& V
}$$
where the diagonal (regular!) maps are obtained by composing the 
rational maps to $V$; and the horizontal map is a blow-up of a
nonsingular variety at a nonsingular center. This center may be chosen
so that it meets with normal crossings the exceptional divisor of the
map to $V$ (cf.~\cite{MR2003c:14016}, expecially the `Furthermore'
section of the main statement, and part~6 of Theorem~0.3.1).

Therefore, at each stage of the factorization we are in the hypotheses
of Proposition~\ref{pre}, with the normal crossing divisor on the base
of the blow-up equal to the Jacobian of the map to $V$; this is
compatible with the prescribed multiplicities $\mu_j$, since the
composition of the differentials is the differential of the
composition. By pasting together the equalities thereby obtained at
each stage, we get that $\chi_Z(\pi^{-1}(U))=\chi_V(U)=[U]$; that is,
the stated equality.
\end{proof}

\section{Proof of the theorem}\label{CSM}

A {\em constructible function\/} on a variety $Z$ is a finite
$\Zbb$-linear combination of characteristic functions of
subvarieties: $\sum_{S\subset Z} m_S \one_S$ with $m_S\in \Zbb$, and
$\one_S(s)=1$ for $s\in S$, $\one_S(s)=0$ otherwise. 

Denote by $F(Z)$ the group of constructible functions on the
variety~$Z$. This assignment is in fact a covariant functor under
proper morphisms, as follows: if $v: Z \to V$ is a proper morphism,
and $S$ is a subvariety (or more generally a constructible set) in
$Z$, define the push-forward of the function $\one_S$ by 
$$\forall p\in V\quad,\quad v_*(\one_S)(p)=\chi(v^{-1}(p)\cap
S)\quad.$$
Here $\chi$ denotes the topological Euler characteristic if the ground
field is $\Cbb$; this definition may be extended to arbitrary
algebraically closed fields of characteristic~0, see
\cite{MR91h:14010}.

In \cite{MR50:13587}, Robert MacPherson proved (over $\Cbb$) that
there exists a natural transformation $c$ from $F$ to homology, such
that if $V$ is a {\em nonsingular\/} variety, then $c(\one_V)$ equals
the total homology Chern class of the tangent bundle of $V$,
$c(TV)\cap [V]$. MacPherson's natural transformation can in fact be
lifted to the Chow group and over any algebraically closed field of
characteristic~0, see Example~19.1.7 in \cite{MR85k:14004} and
\cite{MR91h:14010}.

MacPherson's natural transformation may be used to define a notion
extending to possibly singular varieties the total Chern class of the
tangent bundle. This notion is commonly named
`Chern-Schwartz-MacPherson class', since it can be shown that the
class so obtained agrees with a notion previously defined by
Marie-H\'el\`ene Schwartz, cf.~\cite{MR83h:32011}.

Analogously, we may define a `Chern class' of any constructible set
$S$ in a variety $Z$, {\em as an element of $A_*Z$,\/} by applying
MacPherson's transformation to its characteristic function: that is,
we set 
$$\csm(S):=c(\one_S)\in A_*Z\quad.$$

Now consider the situation of \S\ref{motiv}: $v:Z\to V$ is a proper
birational morphism of nonsingular varieties; the exceptional divisor
has normal crossings, with nonsingular components $E_j$, and the
Jacobian ideal of $v$ is principal, with divisor $\sum_j \mu_j E_j$.

\begin{theorem}\label{strmain}
Let $U$ be a constructible subset of $V$. Then
$$\csm(U)=v_* \sum_{I\subset J}\frac{\csm(E_I^\circ\cap
  v^{-1}(U))}{\prod_{i\in I}(\mu_i+1)}$$
in $(A_*V)_\Qbb$.
\end{theorem}

\begin{proof}
For all $p\in V$,
$$\sum_{I\subset J} \frac{[E_I^\circ\cap v^{-1}(p)]}{\prod_{i\in I}
  [\Pbb^{\mu_i}]}=[p]$$
by Theorem~\ref{basic}. Taking Euler characteristics:
$$\sum_{I\subset J} \frac{\chi(E_I^\circ\cap v^{-1}(p))}{\prod_{i\in I}
  (\mu_i+1)}=1\quad,$$
thus,
$$\sum_{I\subset J} \frac{\chi(E_I^\circ\cap v^{-1}(U)\cap
  v^{-1}(p))}{\prod_{i\in I} (\mu_i+1)}=
\left\{\aligned
0 & \quad p\not\in U\\
1 & \quad p\in U
\endaligned\right.
\quad.$$
By definition of push-forward of constructible functions, this tells
  us that
$$v_*\sum_{I\subset J} \frac{\one_{E_I^\circ\cap v^{-1}(U)}}
{\prod_{i\in I}(\mu_i+1)} =\one_U$$
in $F(V)\otimes\Qbb$. Applying MacPherson's natural transformation
gives the statement.
\end{proof}

Theorem~\ref{main} follows immediately. Indeed, if $V$ and $W$ are
both dominated birationally by~$Z$, and $v: Z \to V$, $w: Z \to W$
have the same Jacobian, we may assume (after resolution of
singularities) that $Z$ is nonsingular, with normal crossing
exceptional divisor $\cup E_j$, and both $v$ and $w$ have Jacobian
with divisor $\sum \mu_j E_j$. Let
$$C=\sum_{I\subset J}\frac{\csm(E_I^\circ)}{\prod_{i\in I}(\mu_i+1)}
\in (A_*Z)_\Qbb\quad;$$
then $\csm(V)=v_*(C)$, $\csm(W)=w_*(C)$ by Theorem~\ref{strmain}, and
these classes agree with $c(TV)\cap [V]$, $c(TW)\cap [W]$
respectively, by the normalization property of
Chern-Schwartz-MacPherson's classes.

Theorem~\ref{strmain} appears to be of independent interest. For
example, let $U$ be a singular variety embedded in a nonsingular
variety $V$, and let $Z$ be the variety obtained from $V$ as $U$
undergoes a resolution of singularities \`a la Hironaka. Then
Theorem~\ref{strmain} gives an explicit expression for the
Chern-Schwartz-MacPherson class of $U$, with rational coefficients, 
in terms of classes of loci in its inverse image in $Z$ (note that
in this case $E_I^\circ\subset v^{-1}(U)$ for $I\ne\emptyset$). 
To our knowledge, this expression is new.

\end{document}